\documentclass{amsproc}%
\usepackage{graphicx}
\usepackage{amscd}
\usepackage{amsmath}
\usepackage{amsfonts}
\usepackage{hyperref}
\usepackage{amssymb}%
\newtheorem{theorem}{Theorem}[section]

\newtheorem{proposition}[theorem]{Proposition}
\theoremstyle{definition}

\theoremstyle{remark}

\numberwithin{equation}{section}
\begin{document}
\title{Robustness of topology of digital images and point clouds}
\author{Peter Saveliev}
\address{Department of Mathematics, Marshall University, One John Marshall Drive,
Huntington, WV 25755}

\begin{abstract}
Such modern applications of topology as digital image analysis and data
analysis have to deal with noise and other uncertainty. In this environment,
the data structures often appear "filtered" into a sequence of cell
complexes. We introduce the homology group of the filtration as the group of
all possible homology classes of all elements of the filtration without
double count. The second step of analysis is to discard the features that
lie outside the user's choice of the acceptable level of noise.
\end{abstract}

\maketitle

\section{Introduction}

Since Poincar\'{e}, homology has been\ used as the main descriptor of the
topology of geometric objects. In the classical context, however, all
homology classes receive equal attention. Meanwhile, applications of
topology in analysis of images and data have to deal with noise and other
uncertainty. This uncertainty appears usually in the form of a real valued
function defined on the topological space. Persistence is a measure of
robustness of the homology classes of the lower level sets of this function 
\cite{ELZ}, \cite{Carlsson}, \cite{CZ09}, \cite{CZ}.

Since it's unknown beforehand what is or is not noise in the dataset, we
need to capture all homology classes including those that may be deemed
noise later. In this paper we introduce an algebraic structure that
contains, without duplication, all these classes. Each of them is associated
with its persistence and can be removed when the acceptable threshold for
noise is set. The last step can be carried out repeatedly in order to find
the best possible threshold. The approach follows the approach to analysis
of digital images presented in \cite{Sav1}.

\section{Backgound}

The topological spaces subject to such analysis are cell complexes. A 
\textit{cell complex} is a combinatorial structure that describes how $k$%
-dimensional cells are attached to each other along $(k-1)$-dimensional
cells. Cell complexes come from the following two main sources.

First, a gray scale image is a real-valued function $f$ defined on a
rectangle. Given a threshold $r$, the lower level set $f^{-1}((-\infty ,r))$
can be thought of as a binary image. Each black pixel of this image is
treated as a square cell in the plane. These 2-dimensional cells are
combined with their edges (1-cells) and vertices (0-cells) and in the $n$%
-dimensional case, the image is decomposed into a combination of $0$-, $1$-,
..., $n$-cubes. This process is called \textit{thresholding}. The result is
a cell complex $K$ for each $r,$ see \cite{KMM}.

Second, a point cloud is a finite set $S$ in some Euclidean space of
dimension $d$. Given a threshold $r$, we deem any two points that lie within 
$r$ from each other as "close". In this case, this pair of points is
connected by an edge. Further, if three points are "close", pairwise, to
each other, we add a face spanned by these points. If there are four, we add
a tetrahedron, and, finally, any $d+1$ "close" points create a $d$-cell. The
process is called the \textit{Vietoris-Rips construction}. The result is a
cell complex $K$ for each $r$ \cite{ELZ}.

Next, we would like to quantify the topology of the cell complex $K.$ It is
done via the \textit{Betti numbersof }$K$: $B_{0}$ is the number of
connected components in $K$; $B_{1}$ is the number of holes or tunnels (1
for letter O or the donut; 2 for letter B and the torus); $B_{2}$ is the
number of voids or cavities (1 for both the sphere and the torus), etc.

The Betti numbers are computed via \textit{homology theory }\cite{Bredon}.
One starts by considering the collection \textit{\ }$C_{k}(K)$ of all formal
linear combinations (over a ring $R$) of $k-$cells in $K,$ called \textit{%
chains}. Combined they form a finitely generated abelian group called the 
\textit{chain complex} $C_{k}(K)$, or collectively $C_{\ast }(K).$ A $k$%
-chain can be recorded as an $N_{k}$-vector, where $N_{k}$ is the total
number of $k$-cells in $K$. The boundary of a $k$-chain is the chain
comprised of all $(k-1)$-faces of its cells taken with appropriate signs.
Then the \textit{boundary operator} $\partial :C_{k}(K)\rightarrow
C_{k-1}(K) $ acts on the chain complex and is represented by a $N_{k}\times
N_{k-1}$ matrix.

From the chain complex $C_{\ast }(K),$ the homology group is constructed by
means of the standard algebraic tools. To capture the topological features
one concentrates on \textit{cycles}, i.e., chains with zero boundary, $%
\partial A=0$. Further, one can verify whether two given $k$-cycles $A$ and $%
B$ are \textit{homologous}: the difference between them is the boundary of a 
$(k+1)$-chain $T:A-B=\partial T$ (such as two meridians of the torus). In
this case, $A$ and $B$ belong to the same \textit{homology class} $H=[A]=[B]$%
. The totality of these equivalence classes in each dimension $k$ is called
the $k$-th \textit{homology group} $H_{k}(K)$ of $K$, collectively $H_{\ast
}(K)$. Then, Betti number $B_{k}$ is the rank of $H_{k}(K)$.

\section{Prior work and outline}

The methods for computing homology groups are well developed. In real-life
applications however both digital images and point clouds may be noisy and
one needs to evaluate the significance of their homology classes. The
approach to this problem has been the following. Instead of using a single
threshold and studying a single cell complex, one considers all thresholds
and all possible cell complexes. Since increasing threshold $r$ enlarges the
corresponding complex, we have a sequence of complexes:%
\begin{equation*}
K^{1}\hookrightarrow K^{2}\hookrightarrow K^{3}\hookrightarrow
K^{4}\hookrightarrow \ldots \ \hookrightarrow K^{s},
\end{equation*}%
where the arrows represent the inclusions: $i^{n,n+1}:K^{n}\hookrightarrow
K^{n+1}.$ Let $i^{nm}:K^{n}\hookrightarrow K^{m},n\leq m,$ also be the
inclusion. This structure $\{K^{n},i^{nm}\}$ is called a \textit{filtration}$%
.$ 

Now, each of these inclusions generates a homomorphism $i_{\ast
}^{nm}:H_{\ast }(K^{n})\rightarrow H_{\ast }(K^{m})$ called the \textit{%
homology map induced by }$i^{nm}.$ As a result, we have a sequence of
homology groups connected by these homomorphisms: 
\begin{equation*}
H_{\ast }(K^{1})\rightarrow H_{\ast }(K^{2})\rightarrow \ldots \ \rightarrow
H_{\ast }(K^{s})\longrightarrow 0.
\end{equation*}%
These homomorphisms record how the homology changes as the complex grows at
each step. For example, a component appears and then merges with another
one, or a hole is formed and then filled. We refer to these events as 
\textit{birth and death} of the corresponding homology classes. 

In order to evaluate the robustness of an element of one of these groups the 
\textit{persistence of a homology class} is defined as the number of steps
it takes for the class to end at $0.$ In other words,%
\begin{equation*}
\text{persistence = death date - birth date.}
\end{equation*}%
The $p$-\textit{persistent homology group} of $K^{i}$ is defined as the
image of $i_{\ast }^{i,i+p}.$ It's what's left from $H_{\ast }(K^{i})$ after 
$p$ steps in the filtration. Now the robustness of the homology classes of
the filtration is evaluated in terms of the set of intervals $[birth,death]$
representing the life-spans, called \textit{barcodes}, of the homology
classes \cite{CZ}.

Our approach is somewhat different. It consists of two steps.

First step: we pool all possible homology classes in all elements of the
filtration together in a single algebraic structure (Sections 4 and 5). The
presence of noise is ignored. The homology group $H_{\ast }(\{K^{n}\})$ of
filtration\textit{\ }$\{K^{n}\}$ captures all homology classes in the whole
filtration -- without double counting. 

Second step: for a given positive integer $p,$ the $p$\textit{-}noise group%
\textit{\ }$N_{\ast }^{p}(\{K^{n}\})$\ is comprised of the homology classes
in $H_{\ast }(\{K^{n}\})$ with the persistence less than $p.$ Next, we
"remove" the noise from the homology group of filtration by using the
quotient (Sections 6 and 7): 
\begin{equation*}
H_{\ast }^{p}(\{K^{n}\})=H_{\ast }(\{K^{n}\})/N_{\ast }^{p}(\{K^{n}\}).
\end{equation*}%
In other words: \textit{if the difference between two homology classes is
deemed noise, they are equivalent}. The second step can be repeated as
needed.

We also discuss the computational aspects of this approach (Section 8) and
multiparameter filtrations (Section 9).

Our approach provides a coarser classification of the homology of
filtrations than the one based on barcodes. The reason is that all homology
classes with long enough life-spans, i.e., high persistence, have equal
place in the homology group $H_{\ast }(\{K^{n}\})$ of the filtration
regardless of the time of birth and death.

\section{Motivation: the homology of a gray scale image}

In this section we will try to understand the meaning of the homology of the
gray scale image in Figure 1. For simplicity we assume that there are only 2
levels of gray in addition to black and white. A visual inspection of the
image suggests that it has three connected components each with a hole.
Therefore, its $0$- and $1$-homology groups shouldhave three generators
each. We now develop an algebraic procedure to arrive at this result.

\begin{figure}[h]
\centering\includegraphics[height=20mm,width=.80\columnwidth]{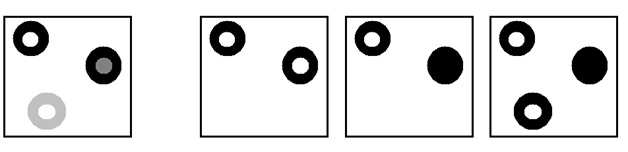} 
\caption{A gray scale image and the corresponding filtration}
\label{fig:image}
\end{figure}

First the image is "thresholded". The lower level sets of the gray scale
function of the image form a filtration: a sequence of three binary images,
i.e. cell complexes: $K^{1}\hookrightarrow K^{2}\hookrightarrow K^{3},$
where the arrows represent the inclusions. Suppose $A_{i},B_{i},C_{i}$ are
the homology classes that represent the components of $K^{i}$ and $%
a_{i},b_{i},c_{i}$ are the holes, clockwise starting at the upper left
corner. The homology groups of these images also form sequences -- one for
each dimension 0 and 1.

Suppose $F_{1},F_{2}$ are\ the two homology maps, i.e., homomorphisms of the
homology groups generated by the inclusions of the complexes, with $F_{3}=0$
included for convenience. These homomorphisms act on the generators, as
follows: 
\begin{eqnarray*}
A_{1} &\rightarrow &A_{2}\rightarrow A_{3}\rightarrow 0,B_{1}\rightarrow
B_{2}\rightarrow B_{3}\rightarrow 0, \\
C_{2} &\rightarrow &C_{3}\rightarrow 0,a_{1}\rightarrow a_{2}\rightarrow
a_{3}\rightarrow 0, \\
b_{1} &\rightarrow &0,c_{3}\rightarrow 0.
\end{eqnarray*}
To avoid double counting, we want to count only the homology classes that
don't reappear in the next homology group. As it turns out, a more
algebraically convenient way to accomplish this is to count only the
homology classes that go to $0$ under these homomorphisms. These classes
form the kernels of $F_{1},F_{2},F_{3}$. Now, we choose the homology group
of the original, gray scale image to be the direct sum of these kernels:%
\begin{equation*}
H_{0}(\{K^{i}\})=<A_{3},B_{3},C_{3}>,\text{ }H_{1}(\{K^{i}%
\})=<b_{1},a_{3},c_{3}>.
\end{equation*}%
Thus the image has three components and three holes, as expected.

\section{Homology groups of filtrations}

In the following sections we provide formal definitions. All cell complexes
are finite.

Suppose we have a one-parameter filtration:%
\begin{equation*}
K^{1}\hookrightarrow K^{2}\hookrightarrow K^{3}\hookrightarrow \ldots \
\hookrightarrow K^{s}.
\end{equation*}%
Here $K^{1},K^{2},\ldots ,K^{s}$ \ are cell complexes and the arrows
represent the inclusions $i^{n,n+1}:K^{n}\hookrightarrow K^{n+1}$ and so do $%
i^{nm}:K^{n}\hookrightarrow K^{m},n\leq m$. We will denote the filtration by 
$\{K^{n},i^{nm}:n,m=1,2,...,s,n\leq m\},$ or simply $\{K^{n}\}.$ Next,
homology generates a "direct system" of groups and homomorphisms:%
\begin{equation*}
H_{\ast }(K^{1})\rightarrow H_{\ast }(K^{2})\rightarrow \ldots \ \rightarrow
H_{\ast }(K^{s})\longrightarrow 0.
\end{equation*}%
We denote this direct system by $\{H_{\ast }(K^{n}),i_{\ast
}^{nm}:n,m=1,2,...,s,n\leq m\},$ or simply $\{H_{\ast }(K^{n})\}.$ The zero
is added in the end for convenience.

Our goal is to define a single structure that captures all homology classes
in the whole filtration without double counting. The rationale is that if $%
x\in H_{\ast }(K^{n}),y\in H_{\ast }(K^{m}),$ $y=i_{\ast }^{nm}(x),$ and
there is no other $x$ satisfying this condition, then $x$ and $y$ may be
thought of as representing the same homology class of the geometric object
behind the filtration. 

The \textit{homology group of filtration }$\{K^{n}\}$ is defined as the
product of the kernels of the inclusions:%
\begin{equation*}
H_{\ast }(\{K^{n}\})=\ker i_{\ast }^{1,2}\oplus \ker \,i_{\ast }^{2,3}\oplus
\ldots \oplus \ker i_{\ast }^{s,s+1}.
\end{equation*}%
Here, from each group we take only the elements that are about to die. Since
each dies only once, there is no double-counting. Since the sequence ends
with $0,$ we know that everyone will die eventually. Hence every homology
class appears once and only once.

These are a few simple facts about this group.

\begin{proposition}
If $i_{\ast}^{n,n+1}$ is an isomorphism for each $n=1,2,...,s-1,$ then
$H_{\ast}(\{K^{n}\})=H_{\ast}(K^{1})$ $.$
\end{proposition}

\begin{proposition}
If $i_{\ast}^{n,n+1}$ is a monomorphism for each $n=1,2,...,s-1,$ then
$H_{\ast}(\{K^{n}\})=H_{\ast}(K^{s}).$
\end{proposition}

\begin{proposition}
Suppose $\{K^{n},i^{nm},n,m=1,2,...,s\}$ and $\{L^{n},j^{nm},n,m=1,2,...,s\}$
are filtrations. Then $H_{\ast}(\{K^{n}\sqcup L^{n}\})=H_{\ast}(\{K^{n}%
\})\oplus H_{\ast}(\{L^{n}\}).$
\end{proposition}

\begin{proposition}
Suppose $\{K^{n},i^{nm},n,m=1,2,...,s\}$ and $\{L^{n},j^{nm},n,m=1,2,...,s\}$
are filtrations and $f:K^{s}\rightarrow L^{s}$ is a cell map. Then the
homology map of the homology groups of these filtrations $f_{\ast}:H_{\ast
}(\{K^{n}\})\rightarrow H_{\ast}(\{L^{n}\})$ is well defined as%
\[
f_{\ast}(x_{1},x_{2},...,x_{s})=(f_{\ast}^{1}(x_{1}),f_{\ast}^{2}%
(x_{2}),...,f_{\ast}^{s}(x_{s})),
\]
where $f^{n}$ is the restriction of $f$ to $K^{n}.$
\end{proposition}

The stability of the homology group of a filtration follows from the
stability of its persistence diagram, i.e., the set of points $%
\{(birth,death)\}\subset \mathbf{R}^{2}$ for the generators of the homology
groups of the filtration, plus the diagonal. It is proven in \cite{CEH} that 
$d_{B}(D(f),D(g))\leq ||f-g||_{\infty },$ where $d_{B}$ is the bottle-neck
distance between the persistence diagrams $D(f),D(g)$ of two filtrations
generated by tame functions $f,g.$ Function $F(x,y)=y-x$ creates an analogue
bottle-neck distance for the set of points $\{persistence\}\subset \mathbf{R}
$ and its stability follows from the continuity of $F$.

\section{Motivation: the high contrast homology of a gray scale image}

To justify our approach to persistence, we observe that some of the features
in the image in Figure 1 are more prominent than others. In particular, some
of the features have lower contrast. These are the holes in the second and
the third rings as well as the third ring itself. By \textit{contrast} of a
lower level set of the gray level function we understand the difference
between the highest gray level adjacent to the set and the lowest gray level
within the set.

An easy computation shows that the homology classes with persistence of 3 or
higher among the generators are: $A_{1},B_{1},a_{1}.$ However, the set of
the classes of high persistence isn't a subgroup of the homology group of
the respective complex. Instead, we look at the classes with \textit{low }%
persistence, i.e., the noise. In particular, the classes in $H_{\ast }(K^{1})
$ of persistence 2 or lower form the kernel of $F_{2}F_{1}$. We now "remove"
this noise from the homology groups of the filtration by considering their
quotients over these kernels. In particular, the 3-persistent homology
groups of the image are:%
\begin{align*}
H_{0}^{3}(\{K^{i}\})& =<A_{1},B_{1}>/0=<A_{1},B_{1}>, \\
H_{1}^{3}(\{K^{i}\})& =<a_{1},b_{1}>/<b_{1}>=<a_{1}>.
\end{align*}%
Observe that the output is identical to the homology of a single complex,
i.e., a binary image, with two components and one hole. The way persistence
is defined ensures that we can never remove a component as noise but keep a
hole in it.

Observe now that the holes in the second and third rings have the same
persistence (contrast) and, therefore, occupy the same position in the
homology group regardless of their birth dates (gray level). Second, if we
shrunk one of these rings, its persistence and, therefore, its place in the
homology group wouldn't change. These observations confirm the fact that the
homology group of the gray scale image, unlike the barcodes, captures only
its topology.

In the case of a Vietoris-Ripps complex, not only the barcode, [birth,
death], but also the persistence, death - birth, of a homology class
contains information about the size of representatives of these classes. For
example, a set of points arranged in a circle will produce a 1-cycle with
twice as large birth, death, and persistence than the same set shrunk by a
factor of 2. However, persistence defined as death/birth will have the
desired property of scale independence. The same result can be achieved by
an appropriate re-parametrizing of the filtration.

\section{Persistent homology groups of filtrations}

In the general context of filtrations the measure of importance of a
homology class is its persistence which is the length of its lifespan within
the direct system of homology of the filtration.

Given filtration $\{K^{n}\},$ we say that \textit{the persistence }$P(x)$ 
\textit{of }$x\in H_{\ast}(K^{n})$\textit{\ is equal to }$p$ if $i_{\ast
}^{n,n+p}(x)=0$ and $i_{\ast}^{n,n+p-1}(x)\neq0.$ Our interest is in the
"robust" homology classes, i.e. the ones with high persistence. However, the
collection of these classes is not a group as it doesn't even contain 0. So
we deal with "noise" first. Given a positive integer $p,$ the $p$\textit{%
-noise (homology) group} $N_{\ast}^{p}(K^{n})$ of $\{K^{n}\}$ is the group
of all elements of $K^{n}$ with persistence less than $p.$

Alternatively, we can define these groups via kernels of the homomorphisms
of the inclusions: $N_{\ast }^{p}(K^{n})=\ker \,i_{\ast }^{n,n+p}.$

\begin{proposition}
$N_{\ast}^{p+1}(K^{n})\subset N_{\ast}^{p}(K^{n}).$
\end{proposition}

Next, we "remove" the noise from the homology group. The $p$\textit{%
-persistent (homology) group} of $K^{n}$ with respect to the filtration $%
\{K^{n}\}$ is defined as 
\begin{equation*}
H_{\ast }^{p}(K^{n})=H_{\ast }(K^{n})/N_{\ast }^{p}(K^{n}).
\end{equation*}%
The point of this definition is that, given a threshold for noise, if the
difference between two homology classes is noise, they should be equivalent.

Next, just as in the case of noise-less analysis, we define a single
structure to capture all (robust) homology classes. Let $p$ be a positive
integer$.$ Suppose $x\in \ker i_{\ast }^{k,k+p}$ and let $y=$ $i_{\ast
}^{k,k+1}(x).$ Then%
\begin{eqnarray*}
i_{\ast }^{k+1,k+1+p}(y) &=&i_{\ast }^{k+1,k+1+p}(i_{\ast }^{k,k+1}(x)) \\
&=&i_{\ast }^{k,k+1+p}(x)=i_{\ast }^{k+p,k+p+1}(i_{\ast }^{k,k+p}(x)) \\
&=&i_{\ast }^{k+p,k+p+1}(0)=0.
\end{eqnarray*}%
Hence $y\in \ker i_{\ast }^{k+1,k+1+p}.$ We have proved that 
\begin{equation*}
i_{\ast }^{k,k+1}(\ker i_{\ast }^{k,k+p})\subset \ker i_{\ast }^{k+1,k+1+p}.
\end{equation*}%
It follows that the homomorphism $i_{\ast }^{k,k+1}:\ker i_{\ast
}^{k,k+p}\rightarrow \ker i_{\ast }^{k+1,k+1+p}$ generated by the inclusion
is well-defined. 

Next, we use these homomorphisms to define the $p$\textit{-noise (homology)
group }$N_{\ast }^{p}(\{K^{n}\})$\textit{\ of filtration} $\{K^{n}\}$ as 
\begin{equation*}
N_{\ast }^{p}(\{K^{n}\})=\ker i_{\ast }^{1,2}\oplus \ldots \oplus \ker
i_{\ast }^{s,s+1}.
\end{equation*}%
Observe that the formula is the same as the one in the definiton of $H_{\ast
}^{p}(\{K^{n}\}).$ Since $i_{\ast }^{k,k+1}:\ker i_{\ast
}^{k,k+p}\rightarrow \ker i_{\ast }^{k+1,k+1+p}$ is a restriction of $%
i_{\ast }^{k,k+1}:H_{\ast }^{p}(K^{k})\rightarrow H_{\ast }^{p}(K^{k+1}),$
each term in the above definition is a subgroup of the corresponding term in
the definition of $H_{\ast }(\{K^{n}\}).$ The proposition below follows.

\begin{proposition}
$N_{\ast }^{p}(\{K^{n}\})\subset H_{\ast }(\{K^{n}\}).$
\end{proposition}

Finally, the\textit{\ }$p$\textit{-persistent (homology) group of filtration 
}$\{K^{n}\}$ is 
\begin{equation*}
H_{\ast }^{p}(\{K^{n}\})=H_{\ast }(\{K^{n}\})/N_{\ast }^{p}(\{K^{n}\}).
\end{equation*}

The results about $H_{\ast }^{p}(\{K^{n}\})$ analogous to the ones about $%
H_{\ast }(\{K^{n}\})$ in Section 5 hold.

\section{Computational aspects}

For 2-dimensional gray scale images, this approach to homology and
persistence has been used in an image analysis program. The algorithm
described in \cite{Sav1} has complexity of $O(n^{2}),$ where $n$ is the
number of pixels in the image.

For the general case, the analysis algorithm may be outlined as follows:

\begin{enumerate}
\item The input is a filtration.

\item The homology groups of its members and the homomorphisms induced by
inclusions are computed.

\item The homology group of the filtration is computed.

\item The persistence of all elements of the homology groups is computed.

\item The user sets a threshold $p$ for persistence and the $p$-noise group
of the filtration is computed.

\item The $p$-persistent homology group of the filtration is computed and
given as output.
\end{enumerate}

If the user changes the threshold, the last step is repeated as necessary
without repeating the rest.

The algorithm above computes the homology group of filtration, as defined,
incrementally. This may be both a disadvantage and an advantage. In
comparison, the persistence complex \cite{CZ} also contains information
about all homology classes of the filtration but its computation does not
require computing the homology of each complex of the filtration. Meanwhile,
the above algorithm may have to compute the same homology over and over if
consecutive complexes are identical. Hence, the algorithm has a disadvantage
in terms of processing time. On the other hand, the incremental nature of
the algorithm makes its use\ of memory independent from the length of the
filtration. Another advantage is that multi-parameter filtrations are dealt
with in the exact same manner (see next section).

The inefficiency of the above algorithm can be addressed with a proper
algebraic tool. This tool is the mapping cone \cite{Weibel}. Suppose, for
simplicity, that our filtration has only two elements: $i:K^{1}%
\hookrightarrow K^{2}.$ The mapping cone is, in a sense, a combination of
the kernel and the cokernel of $i_{\ast }$. It captures the difference
between $K^{1}$ and $K^{2}$ on the chain level: everything in $C_{\ast
}(K^{1})$ is killed unless it also appears in $C_{\ast }(K^{2})$ under $%
i_{\ast }$. Then the algorithm is to construct the homology group from the
chain complexes $C_{\ast }(K^{1}),C_{\ast }(K^{2})$ of the elements of the
filtration and the chain map $i_{\ast }:C_{\ast }(K^{1})\rightarrow C_{\ast
}(K^{2}).$

\section{Multiparameter filtrations}

Multiparameter filtrations come from the same main sources as one-parameter
filtrations. First, color images are thresholded according to their three
color channels. Second, point clouds are thresholded by the closeness of
their points and, for example, the density of hte points. 

Let limit our attention to the two-parameter case. A (finite) two-parameter
filtrations $\{K^{nm}\}$ is a table of complexes connected by inclusions 
\begin{equation*}
i(n,m,n+p,m+q):K^{nm}\rightarrow K^{n+p,m+q},p,q\geq 0,
\end{equation*}
These inclusions generate homomorphisms 
\begin{equation*}
i_{\ast }(n,m,n+q,m+p):H_{\ast }(K^{nm})\rightarrow H_{\ast }(K^{n+q,m+p}),
\end{equation*}%
with 0s added in the end of each row and each column. Define the \textit{%
homology group of the filtration }$\{K^{nm}\}$ as%
\begin{eqnarray*}
H_{\ast }(\{K^{nm}\}) = {\displaystyle\bigoplus\limits_{n}}\ker i_{\ast }(n,m,n+1,m)\cap \ker i_{\ast
}(n,m,n,m+1).
\end{eqnarray*}%
The analogues of the results in Section 5 hold.

There are many ways to define persistence in the multiparameter setting. For
example, we can evaluate the robustness of a homology class $x\in H_{\ast
}(K^{nm})$ in terms of the pairs $(p,q)$ of positive integers satisfying%
\begin{equation*}
i_{\ast }(n,m,n+p,m)(x)=0\text{ and }i_{\ast }(n,m,n,m+q)(x)=0.
\end{equation*}%
Next, just as in Section 7, we restrict the homomorphisms generated by the
inclusions to the homology classes of low persistence:%
\begin{eqnarray*}
i_{\ast }(n,m,n+1,m) &:& \\
\ker i_{\ast }(n,m,n+p,m) &\rightarrow &\ker i_{\ast }(n+1,m,n+1+p,m), \\
i_{\ast }(n,m,n,m+1) &:& \\
\ker i_{\ast }(n,m,n,m+q) &\rightarrow &\ker i_{\ast }(n+1,m,n,m+1+q).
\end{eqnarray*}%
Then the $(p,q)$-\textit{noise group of} $K^{nm}$ is defined via these
homomorphisms:%
\begin{eqnarray*}
N_{\ast }^{pq}(\{K^{nm}\})  = {\displaystyle\bigoplus\limits_{n}}\ker i_{\ast }(n,m,n+1,m)\cap \ker i_{\ast
}(n,m,n,m+1).
\end{eqnarray*}%
Finally, the $(p,q)$\textit{-persistent (homology) group of filtration }$%
\{K^{nm}\}$ is defined as 
\begin{equation*}
H_{\ast }^{pq}(\{K^{n}\})=H_{\ast }(\{K^{nm}\})/N_{\ast }^{pq}(\{K^{nm}\}).
\end{equation*}

The results about $H_{\ast }^{pq}(\{K^{nm}\})$ analogous to the ones about $%
H_{\ast }^{p}(\{K^{n}\})$ in Section 7 hold.


\href{http://inperc.com/wiki/index.php?title=Homology_of_filtrations}{Robustness of topology of digital images and point clouds}

\end{document}